\documentclass[11pt,reqno]{amsart}

\usepackage[a4paper,left=35mm,right=35mm,top=30mm,bottom=30mm,marginpar=25mm]{geometry}

\usepackage{amsfonts}
\usepackage{amsmath}
\usepackage{amssymb}
\usepackage{amsthm}
\usepackage[latin1]{inputenc}
\usepackage{eurosym}
\usepackage[dvips]{graphics}
\usepackage{graphicx}
\usepackage{epsfig}
\usepackage{hyperref}
\usepackage{dsfont}
\usepackage{color}
\usepackage[dvipsnames]{xcolor}
\usepackage{booktabs, multirow} 
\usepackage{soul}
\usepackage{changepage,threeparttable} 
\usepackage{graphicx,subcaption}
\usepackage[displaymath,mathlines]{lineno}

\allowdisplaybreaks

\usepackage{hyperref}

\usepackage{mathtools}

\usepackage{ifthen}
\graphicspath{./figures/}
\makeindex

\newcommand{\la}{{\langle}}
\newcommand{\ra}{{\rangle}}

\newcommand{\Pp}{\mathbb{P}}

\newcommand{\bK}{\mathbf{K}}
\newcommand{\bS}{\mathbf{S}}

\numberwithin{equation}{section}

\newtheoremstyle{thmlemcorr}{10pt}{10pt}{\itshape}{}{\bfseries}{.}{10pt}{{\thmname{#1}\thmnumber{
			#2}\thmnote{ (#3)}}}
\newtheoremstyle{thmlemcorr*}{10pt}{10pt}{\itshape}{}{\bfseries}{.}\newline{{\thmname{#1}\thmnumber{
			#2}\thmnote{ (#3)}}}
\newtheoremstyle{defi}{10pt}{10pt}{\itshape}{}{\bfseries}{.}{10pt}{{\thmname{#1}\thmnumber{
			#2}\thmnote{ (#3)}}}
\newtheoremstyle{remexample}{10pt}{10pt}{}{}{\bfseries}{.}{10pt}{{\thmname{#1}\thmnumber{
			#2}\thmnote{ (#3)}}}
\newtheoremstyle{ass}{10pt}{10pt}{}{}{\bfseries}{.}{10pt}{{\thmname{#1}\thmnumber{
			A#2}\thmnote{ (#3)}}}

\theoremstyle{thmlemcorr}
\newtheorem{theorem}{Theorem}
\numberwithin{theorem}{section}

\newtheorem{proposition}[theorem]{Proposition}

\theoremstyle{thmlemcorr*}
\newtheorem{theorem*}{Theorem}
\newtheorem{lemma*}[theorem]{Lemma}
\newtheorem{corollary*}[theorem]{Corollary}
\newtheorem{proposition*}[theorem]{Proposition}
\newtheorem{problem*}[theorem]{Problem}
\newtheorem{conjecture*}[theorem]{Conjecture}

\theoremstyle{defi}

\theoremstyle{remexample}
\newtheorem{remark}[theorem]{Remark}

\theoremstyle{ass}

\begin{document}
	
	\title[Splitting methods for a class of non-potential mean field games]{Splitting methods for a class of non-potential mean field games}
	
	
	\author{Siting Liu}
	\address[S. L.]{Department of Mathematics, UCLA}
	\email{siting6@math.ucla.edu}
	
	
	
	\author{Levon Nurbekyan}
	\address[L. N.]{Department of Mathematics, UCLA}
	\email{lnurbek@math.ucla.edu}
	
	
	\keywords{}
	\subjclass[2010]{Primary: 35Q91, 35Q93, 65M70 Secondary: 35A15, 93A15}
	

	\thanks{This work was supported by AFOSR MURI FA9550-18-1-0502, AFOSR Grant No. FA9550-18-1-0167, ONR Grant No. N00014-18-1-2527, ONR N00014-20-1-2093.} 
	\date{\today}
	
\begin{abstract}
We extend the methods from \cite{nursaude18,liu2020computational} to a class of \textit{non-potential} mean-field game (MFG) systems with mixed couplings. Up to now, splitting methods have been applied to \textit{potential} MFG systems that can be cast as convex-concave saddle-point problems. Here, we show that a class of non-potential MFG can be cast as primal-dual pairs of monotone inclusions and solved via extensions of convex optimization algorithms such as the primal-dual hybrid gradient (PDHG) algorithm. A critical feature of our approach is in considering dual variables of nonlocal couplings in \textit{Fourier} or \textit{feature spaces}.
\end{abstract}
	
	\maketitle
	

\section{Introduction}
	
Our goal is to develop computational methods for the mean-field games (MFG) systems of the form
\begin{equation}\label{eq:main}
\begin{cases}
-\phi_t+H(t,x,\nabla\phi)=f_0(t,x,\rho(x,t))+f_1\left(t,x,\int_\Omega K(x,y) \rho(y,t)dy\right)\\
\rho_t-\nabla \cdot \left(\rho \nabla_p H(t,x,\nabla\phi(x,t))\right)=0\\
\rho(x,0)=\rho_0(x),~ \phi(x,1)=g_0(x,\rho(x,1))+g_1 \left(x,\int_\Omega S(x,y)\rho(y,1)dy\right)
\end{cases}
\end{equation}
This system characterizes Nash equilibria of a differential game with a continuum of agents. For a detailed introduction and description of these models we refer to seminal papers \cite{LasryLions06a,LasryLions06b,LasryLions2007,HCM06,HCM07}, manuscripts \cite{gueantlasrylions11,CardaNotes,gomsaude'14}, and references therein.

In \eqref{eq:main}, $\phi$ is the value function of a generic agent, $\rho$ is the distribution of the agents in the state-space, $H$ is the Hamiltonian of a single agent, and $f_0,f_1,g_0,g_1$ are the terms that model interactions between a single agent and the population. These interactions can be either local, $f_0,g_0$ terms, or nonlocal, $f_1,g_1$ terms. In the latter case, one needs to assemble information across the whole population using interaction kernels $K,S$. More specifically, $K(x,y),S(x,y)$ signify how agents located at $y$ affect the decision-making of an agent at $x$. Finally, $\Omega \subset \mathbb{R}^d$ is a smooth domain.

Computational methods for \eqref{eq:main} can be roughly divided into two groups. The first group of methods applies to a specific class of MFG systems that are called \textit{potential} and can be cast as convex-concave saddle-point problems. In this case, \eqref{eq:main} can be efficiently solved via splitting methods such as alternating direction method of multipliers (ADMM) \cite{bencar'15,bencarsan'17} and primal-dual hybrid gradient (PDHG) algorithm \cite{silva18,silva19}. These methods are mostly applicable to systems with only local couplings because only a limited class of problems with nonlocal ones are potential: see Section 1.2 in \cite{hadi17}.

The second group of methods are general purpose and do not rely on a specific structure of \eqref{eq:main}. We refer to \cite{achdou10,achdou13,achdou13b} for finite-difference, \cite{silva12,silva14,silva15,carlini18} for semi-Lagrangian, \cite{gomes17,saude18} for monotone flow, and \cite{hadi17,hadi17b,hadi19} for game-theoretic learning methods.

All of the methods above, work directly with discretizations of interaction terms in the state-space. Therefore calculations of nonlocal terms require storing discretizations of $K,S$ and assembling $\int_\Omega K(x,y) \rho(y,t) dy,~\int_\Omega S(x,y) \rho(y,t) dy$ via matrix products across the whole grid. Hence, this procedure is prone to high memory and computational costs, especially on a fine grid. In \cite{nursaude18,liu2020computational}, the authors remedy this issue by passing to Fourier coordinates in the nonlocal terms. Relying on approximations of $K,S$ in a suitable basis, they approximate nonlocal terms with a relatively small number of parameters independent of the grid-size. Additionally connections are drawn with kernel methods in machine-learning \cite[Section 4]{liu2020computational}.

Systems considered in \cite{nursaude18,liu2020computational} have only nonlocal interactions and are potential. Here, we extend these results to systems of the form \eqref{eq:main} that are non-potential in general and contain both local and nonlocal interactions. Our method relies on a monotone-inclusion formulation of \eqref{eq:main} where inputs from different interaction terms are split via dual variables. A critical feature of the method is that dual variables corresponding to nonlocal terms are set up in Fourier spaces.

We list a number of advantages of our method. Firstly, the Fourier approach yields a dimension reduction: see \cite[Section 3.1]{liu2020computational} for a detailed discussion. Secondly, any number of local and nonlocal interactions can be added to the system bearing minimal and straightforward changes on the algorithm. Furthermore, the algorithm is highly modular and parallelizable. Indeed, the updates of dual variables corresponding to different interactions are decoupled. Finally, the monotone inclusion formulation readily provides the convergence guarantees for the algorithm.

Here, we do not concentrate on theoretical aspects of \eqref{eq:main} and our derivations are mostly formal. We refer to \cite{alpar18} and references therein for rigorous treatment of systems similar to \eqref{eq:main}. The paper is organized as follows. In Section \ref{sec:mon_inclusion}, we present our approach and derive the monotone-inclusion formulation of \eqref{eq:main}. In Section \ref{sec:algo}, we propose a primal-dual algorithm based on this formulation. Furthermore, in Section \ref{sec:kernels}, we consider a concrete class of non-potential models with density constraints. Next, in Section \ref{sec:numerics}, we provide numerical examples. Finally, Appendix contains some of the formal derivations.

\section{MFG via monotone inclusions}\label{sec:mon_inclusion}

We solve \eqref{eq:main} in two steps. Firstly, we approximate \eqref{eq:main} by a lower-dimensional system via orthogonal projections of nonlocal terms. Next, we formulate the lower-dimensional system as a monotone inclusion problem. 

Assume that $\{\zeta_i\}_{i=1}^r\subset C^2(\Omega)$ is an orthonormal system with respect to the $L^2(\Omega)$ inner product. Then the system of functions $\{\zeta_i \otimes \zeta_j\}_{i,j=1}^r$ is also orthonormal, where $\zeta_i \otimes \zeta_j (x,y)=\zeta_i(x) \zeta_j(y)$. Furthermore, denote by $\Pp_r$ and $\Pp_{r,r}$ the orthogonal projection operators in $L^2(\Omega)$ and $L^2(\Omega^2)$ onto $\mathrm{span}\{\zeta_i\}_{i=1}^r$ and $\mathrm{span}\{\zeta_i \otimes \zeta_j \}_{i,j=1}^r$, respectively. Now consider the following approximation of \eqref{eq:main}
\begin{equation}\label{eq:main_r}
\begin{cases}
-\phi_t+H(t,x,\nabla\phi)=f_0(t,x,\rho(x,t))+\Pp_r \left( f_1\left(t,\cdot,\int_\Omega K_r(\cdot,y) \rho(y,t)dy\right) \right)(x),\\
\rho_t-\nabla \cdot \left(\rho \nabla_p H(t,x,\nabla \phi(x,t))\right)=0,\\
\rho(x,0)=\rho_0(x),~ \phi(x,1)=g_0(x,\rho(x,1))+\Pp_r \left(g_1 \left(\cdot,\int_\Omega S_r(\cdot,y)\rho(y,1)dy\right)\right)(x),
\end{cases}
\end{equation}
where $K_r=\Pp_{r,r}(K),~S_r=\Pp_{r,r}(S)$. For smooth $K,S,f_1,g_1$ and a suitable choice of $\{\zeta_i\}_{i=1}^r$, solutions of \eqref{eq:main_r} approximate those of \eqref{eq:main}. Furthermore, we solve \eqref{eq:main_r} by the \textit{coefficients method} proposed in \cite{nurbekyan18,nursaude18,liu2020computational}. The key observation is that for any $\rho$ we a priori have that
\begin{equation}\label{eq:P_r_f,g}
\begin{split}
	\Pp_r \left(f_1\left(t,\cdot,\int_\Omega K_r(\cdot,y) \rho(y,t)dy\right)\right)(x)=&\sum_{i=1}^{r} a_i(t) \zeta_i(x)\\
	\Pp_r \left(g_1\left(\cdot,\int_\Omega S_r(\cdot,y) \rho(y,1)dy\right)\right)(x)=&\sum_{i=1}^{r} b_i \zeta_i(x),
\end{split}
\end{equation}
where
\begin{equation*}
\begin{split}
a_i(t)=&\int_\Omega f_1\left(t,x,\int_\Omega K_r(x,y) \rho(y,t)dy\right) \zeta_i(x)dx,\\
b_i=&\int_\Omega g_1\left(x,\int_\Omega S_r(x,y) \rho(y,1)dy\right) \zeta_i(x)dx.
\end{split}
\end{equation*}
Therefore, introducing variables
\begin{equation}\label{eq:alpha_beta}
	\alpha(x,t)=f_0(t,x,\rho(x,t)), \quad \beta(x)=g_0(x,\rho(x,1)),\quad m(x,t)=\rho \nabla_p H(t,x,\nabla \phi),
\end{equation}
we obtain that \eqref{eq:main_r} is equivalent to
\begin{equation}\label{eq:main_r_alpha}
	\begin{cases}
	-\phi_t+H(t,x,\nabla\phi)=\alpha(x,t)+\sum_{i=1}^{r} a_i(t) \zeta_i(x),\\
	m(x,t)=-\rho \nabla_p H(t,x,\nabla \phi),\\
	\rho(x,0)=\rho_0(x),~ \phi(x,1)=\beta(x)+\sum_{i=1}^{r} b_i \zeta_i(x),
	\end{cases}
\end{equation}
supplemented with compatibility conditions
\begin{equation}\label{eq:consistency}
\begin{cases}
a_i(t)=\int_\Omega f_1\left(t,x,\int_\Omega \Pp_{r,r}(K)(x,y) \rho(y,t)dy\right) \zeta_i(x)dx,~\forall i\\
b_i=\int_\Omega g_1\left(x,\int_\Omega \Pp_{r,r}(S)(x,y) \rho(y,1)dy\right) \zeta_i(x)dx,~\forall i\\
\alpha(x,t)=f_0(t,x,\rho(x,t))\\
\beta(x)=g_0(x,\rho(x,1))\\
\rho_t+\nabla \cdot m=0
\end{cases}
\end{equation}
Thus, our goal is to formulate \eqref{eq:main_r_alpha}-\eqref{eq:consistency} as a monotone inclusion problem. For that, we start by casting \eqref{eq:main_r_alpha} as a convex duality relation between variables $\left(a, b, \alpha, \beta, \phi\right)$ and $(\rho,m)$. We omit the domains in the notation of Lebesgue and Sobolev spaces when there is no ambiguity. Recall that
\begin{equation*}
	L(t,x,v)=\sup_p -v\cdot p - H(t,x,p)
\end{equation*}
is the convex dual of $H$.

\begin{proposition}\label{prp:Psi}
	For $(a,b,\alpha,\beta,\phi)\in L^2_t \times l^2 \times L^2_{x,t}\times L^2_x \times H^1_{x,t}$ define
	\begin{equation}\label{eq:Psi}
	\begin{split}
	\Psi(a,b,\alpha,\beta,\phi)=&\inf_{\rho,m} \Xi(\rho,m)+\int_{0}^1 \int_\Omega \left(\alpha(x,t)+\sum_{i=1}^r a_i(t)\zeta_i(x)\right)\rho(x,t) dxdt\\
	&+\int_\Omega \left(\beta(x)+\sum_{i=1}^r b_i \zeta_i(x)\right)\rho(x,1)dx\\
	&+\int_0^1 \int_\Omega \phi_t \rho +\nabla \phi \cdot mdxdt-\int_{\Omega} \phi(x,1) \rho(x,1)dx\\
	&+\int_{\Omega} \phi(x,0) \rho(x,0)dx,
	\end{split}
	\end{equation}
	where
	\begin{equation}\label{eq:Xi}
		\Xi(\rho,m)=\int_0^1 \int_\Omega \rho L\left(t,x,\frac{m}{\rho}\right)dxdt+\mathbf{1}_{\rho(x,0)=\rho_0(x)}+\mathbf{1}_{\rho\geq 0}+\mathbf{1}_{m<<\rho}.
	\end{equation}
	Then $(\phi,\rho)$ satisfy \eqref{eq:main_r_alpha} if and only if $(\rho,m)$ is a solution of the optimization problem in \eqref{eq:Psi}. Furthermore,
	\begin{equation}\label{eq:Psi_Xi*}
		\Psi(a,b,\alpha,\beta,\phi)=-\Xi^*\left(C(a,b,\alpha,\beta,\phi)\right),
	\end{equation}
	where
	\begin{equation}\label{eq:C}
		C(a,b,\alpha,\beta,\phi)=\begin{pmatrix}
		-\sum_i a_i(t) \zeta_i(x)-\alpha(x,t)-\phi_t\\
		-\nabla \phi\\
		-\phi(x,0)\\
		-\sum_i b_i \zeta_i(x)-\beta(x)+\phi(x,1)
		\end{pmatrix},
	\end{equation}
	and $\Xi^*$ is the convex dual of $\Xi$; that is,
	\begin{equation*}
	\begin{split}
	\Xi^*(\hat{\rho},\hat{m},\hat{\rho}(\cdot,0),\hat{\rho}(\cdot,1))=&\sup_{\rho,m} \int_{0}^1 \int_\Omega \hat{\rho} \rho+\hat{m}\cdot m dx dt+\int_\Omega \hat{\rho}(x,0)\rho(x,0)dx\\
	&+\int_\Omega \hat{\rho}(x,1)\rho(x,1)dx-\Xi(\rho,m).
	\end{split}
	\end{equation*}
	In particular, $(\rho,m)$ is a solution of the optimization problem \eqref{eq:Psi} if and only if
	\begin{equation}\label{eq:Calpha_rho}
		C(a,b,\alpha,\beta,\phi)\in \partial \Xi(\rho,m),\quad \mbox{or} \quad (\rho,m)\in \partial \Xi^*\left(C(a,b,\alpha,\beta,\phi)\right)
	\end{equation}
\end{proposition}
\begin{proof}
See Appendix.
\end{proof}

Next step is to find a map $M$ such that the relation \eqref{eq:consistency} can be written as $-\left(C^*(\rho,m,\rho(\cdot,0),\rho(\cdot,1))\right) \in M(a,b,\alpha,\beta,\phi) $, where $C^*$ is the adjoint operator of $C$, and $M$ is a maximally monotone map. We have that
\begin{equation}\label{eq:K_S_expansion}
\begin{split}
K_r(x,y)=&\sum_{p,q=1}^r k_{pq} \zeta_p(x)\zeta_q(y),\quad k_{pq}=\int_{\Omega^2} K(x,y)\zeta_p(x)\zeta_q(y)dxdy,\\
S_r(x,y)=&\sum_{p,q=1}^r s_{pq} \zeta_p(x)\zeta_q(y),\quad s_{pq}=\int_{\Omega^2} S(x,y)\zeta_p(x)\zeta_q(y)dxdy,
\end{split}
\end{equation}
and we denote by $\bK=(k_{pq}),~\bS=(s_{pq})$. Without loss of generality, we assume that $\bK,~\bS$ are invertible.

Additionally, assume that $f_0(t,x,\cdot)$, $g_0(x,\cdot)$, $f_1(t,x,\cdot)$, $g_1(x,\cdot)$ are increasing. This assumption means that agents are crowd averse that leads to a well posed system \eqref{eq:main} \cite{LasryLions2007}. Furthermore, denote by
\begin{equation}\label{eq:U}
\begin{split}
U_0(\rho)=\int_0^1 \int_\Omega F_0(t,x,\rho(x,t))dx dt,\quad& V_0(\mu)=\int_\Omega G_0(x,\mu(x)) dx \\
U_1(c)=\int_0^1 \int_\Omega F_1\left(t,x,\sum_p c_p(t)\zeta_p(x)\right)dx dt,\quad& V_1(w)=\int_\Omega G_1\left(x,\sum_p w_p \zeta_p(x)\right)dx,
\end{split}
\end{equation}
where $\partial_z \Diamond_i(t,x,z)=\Diamond_i$ for $\Diamond \in \left\{f,g\right\}$ and $i\in \{0,1\}$. Then we have that $U_i,V_i$ are convex, and we can consider their dual functions
\begin{equation}\label{eq:U*}
\begin{split}
U_1^*(a)=\sup_c \int_0^1 \sum_p a_p(t) c_p(t) dt-U_1(c),\quad V_1^*(b)=\sup_w b\cdot w-V_1(w),\\
U_0^*(\alpha)=\sup_\rho \int_0^1 \int_{\Omega} \alpha \rho dxdt-U_0(\rho),\quad V_0^*(\beta)=\sup_\mu \int_{\Omega} \beta \mu dx-V_0(\mu).
\end{split}
\end{equation}
\begin{proposition}\label{prp:M}
	Assume that $C$ is defined as in \eqref{eq:C}. Then we have that
	\begin{equation}\label{eq:C*}
		C^*(\rho,m,\rho(\cdot,0),\rho(\cdot,1))=\begin{pmatrix}
		(-\int_\Omega \rho(x,t)\zeta_i(x)dx)_i\\
		(-\int_\Omega \rho(x,1)\zeta_i(x)dx)_i\\
		-\rho(x,t)\\
		-\rho(x,1)\\
		- \mathcal{L}^{-1}(\rho_t+\nabla \cdot m,0,0)\\
		\end{pmatrix},
	\end{equation}
	where $\mathcal{L}=\left(\Delta_{t,x},\left(\operatorname{Id}-\partial_t\right)\lfloor_{\Omega\times\{0\}},\left(\operatorname{Id}+\partial_t\right)\lfloor_{\Omega\times\{1\}}\right)$. Furthermore, \eqref{eq:consistency} is equivalent to
	\begin{equation}\label{eq:C*rho_alpha}
		-\left(C^*(\rho,m,\rho(\cdot,0),\rho(\cdot,1))\right) \in M(a,b,\alpha,\beta,\phi), 
	\end{equation}
	where
	\begin{equation}\label{eq:M}
			M(a,b,\alpha,\beta,\phi)=\begin{pmatrix}
		\bK^{-1}~\partial_a U_1^*(a)\\
		\bS^{-1}~\partial_b V_1^*(b)\\
		\partial_\alpha U_0^*(\alpha)\\
		\partial_\beta V_0^*(\beta)\\
		0
		\end{pmatrix}
	\end{equation}
\end{proposition}
\begin{proof}
See Appendix.
\end{proof}
\begin{remark}
The inverse Laplacian operator appears in $C^*$ because we consider $\phi$ as an element of $H^1$ space rather than $L^2$. As an inner product in $H^1$ we set
\begin{equation*}
\la \phi,h \ra_{H^1}=\int_\Omega \phi(x,0) h(x,0) dx+\int_\Omega \phi(x,1) h(x,1) dx+\int_0^1 \int_\Omega \nabla_{t,x} \phi \cdot \nabla_{t,x} h dx dt
\end{equation*}
As pointed out in \cite{matt19,mattflavien}, the choice of spaces is crucial for grid-independent convergence of primal-dual algorithms. We come back to this point below when we discuss the algorithms.
\end{remark}

Combining Propositions \ref{prp:Psi}, \ref{prp:M} we obtain the following theorem.

\begin{theorem}\label{thm:inclusion}
The pair of systems \eqref{eq:main_r_alpha}-\eqref{eq:consistency}, and so \eqref{eq:main_r}, can be written as a primal-dual pair of inclusions
\begin{equation}\label{eq:inclusion_ab_phi}
\begin{split}
0\in M(a,b,\alpha,\beta,\phi)+C^*(N(C(a,b,\alpha,\beta,\phi))) \quad & \mbox{(P)}\\
\begin{cases}
(\rho,m,\rho(\cdot,0),\rho(\cdot,1)) \in N(C(a,b,\alpha,\beta,\phi))\\
-C^*(\rho,m,\rho(\cdot,0),\rho(\cdot,1)) \in M(a,b,\alpha,\beta,\phi)
\end{cases}
\quad & \mbox{(D)}
\end{split}
\end{equation}
where $N=\partial \Xi^*$. Furthermore, if $c \mapsto \partial_c U_1(\bK c)$ and $w\mapsto \partial_w V_1(\bS w)$ are maximally monotone, then $M$ is maximally monotone, and \eqref{eq:inclusion_ab_phi} is a primal-dual pair of monotone inclusions.
\end{theorem}
\begin{proof}
See Appendix.
\end{proof}
\begin{remark}
The monotonicity of $c \mapsto \partial_c U_1(\bK c)$ and $w\mapsto \partial_w V_1(\bS w)$ yields that the mean-field coupling in \eqref{eq:main_r} satisfies the Lasry-Lions monotonicity condition \cite[Theorems 2.4, 2.5]{LasryLions2007}, and hence \eqref{eq:main_r} is well-posed.
\end{remark}

\section{A monotone primal-dual algorithm}\label{sec:algo}

We apply the monotone primal-dual algorithm in \cite{vu13} to solve \eqref{eq:inclusion_ab_phi} (and thus \eqref{eq:main_r}). We start by an abstract discussion of the algorithm. Following \cite{vu13}, assume that $\mathcal{H}, \mathcal{G}$ are real Hilbert spaces, $M:\mathcal{H} \to 2^{\mathcal{H}}$, $N:\mathcal{G} \to 2^{\mathcal{G}}$ are maximally monotone operators, and $C:\mathcal{H} \to \mathcal{G}$ is a nonzero bounded linear operator. Furthermore, consider the following pair of monotone inclusion problems
\begin{equation}\label{eq:gen_mon}
\begin{split}
\mbox{find}~s~\mbox{s.t.}~0\in M s+C^*(N(Cs))\quad &\mbox{(P)}\\
\mbox{find}~q~\mbox{s.t.}~q\in N(Cs),~-C^*q \in Ms,~\mbox{for some}~s\quad &\mbox{(D)}
\end{split}
\end{equation}
When $M=\partial f$, $N=\partial g$ \eqref{eq:gen_mon} reduces to a convex-concave saddle-point problem
\begin{equation*}
	\inf_s f(s)+g(Cs)=  \inf_s\sup_q f(s)+\la C s, q \ra - g^*(q)
\end{equation*}
Accordingly, one can solve \eqref{eq:gen_mon} by a monotone-inclusion version of the celebrated primal-dual hybrid gradient (PDHG) method \cite{champock11,champock16}. In its simplest form, the algorithm in \cite{vu13} reads as follows
\begin{equation}\label{eq:monotone_pdhg}
	\begin{cases}
	s^{n+1}=J_{\tau_s M}(s^n-\tau_s C^*q^n)\\
	\tilde{s}^{n+1}=2s^{n+1}-s^n\\
	q^{n+1}=J_{\tau_q N^{-1}}(q^n+\tau_q C\tilde{s}^{n+1}),
	\end{cases}
\end{equation}
where $J_{\tau F}=(I+\tau F)^{-1}$ is the resolvent operator, and $\tau_s,\tau_q>0$ are such that $ \tau_s \tau_q \|C\|^2 < 1$. Note that when $M=\partial f$, $N=\partial g$ \eqref{eq:monotone_pdhg} reduces to the standard PDHG \cite{champock11,champock16}.

\subsection{A primal-dual algorithm}

Applying \eqref{eq:monotone_pdhg} to \eqref{eq:inclusion_ab_phi} we obtain the following algorithm to solve the MFG system \eqref{eq:main_r}: 
\begin{equation}\label{eq:algorithm}
	\begin{cases}
	(a^{n+1},b^{n+1},\alpha^{n+1},\beta^{n+1},\phi^{n+1})\\
	=J_{\tau M}\bigg((a^{n},b^{n},\alpha^{n},\beta^{n},\phi^{n})-\tau C^*(\rho^n,m^n,\rho^n(\cdot,0),\rho^n(\cdot,1))\bigg)\\
	(\tilde{a}^{n+1},\tilde{b}^{n+1},\tilde{\alpha}^{n+1},\tilde{\beta}^{n+1},\tilde{\phi}^{n+1})\\
	=2(a^{n+1},b^{n+1},\alpha^{n+1},\beta^{n+1},\phi^{n+1})-(a^{n},b^{n},\alpha^{n},\beta^{n},\phi^{n})\\
	(\rho^{n+1},m^{n+1},\rho^{n+1}(\cdot,0),\rho^{n+1}(\cdot,1))\\
	=J_{\sigma \partial \Xi}\bigg((\rho^{n},m^{n},\rho^{n}(\cdot,0),\rho^{n}(\cdot,1))+\sigma C(\tilde{a}^{n+1},\tilde{b}^{n+1},\tilde{\alpha}^{n+1},\tilde{\beta}^{n+1},\tilde{\phi}^{n+1}) \bigg)
	\end{cases}
\end{equation}

\begin{remark}
The time-steps in \eqref{eq:algorithm} must satisfy the condition $\tau \sigma \|C\|^2< 1$. Note that
\begin{equation*}
\begin{split}
\left|\la C(a,b,\alpha,\beta,\phi),(\rho,m,\rho(\cdot,0 ),\rho(\cdot,1) ) \ra \right| \leq &\left(\|a\|_{L^2_t}+\|b\|_2+\|\alpha\|_{L^2_{x,t}}+\|\beta\|_{L^2_x}+\|\phi\|_{H^1}\right)\\
& \left(\|\rho\|_{L^2_{x,t}}+\|m\|_{L^2_{x,t}}+\|\rho(\cdot,0)\|_{L^2_x}+\|\rho(\cdot,1)\|_{L^2_x}\right)
\end{split}
\end{equation*}
Therefore, $\|C\|$ is finite, and independent of the grid-size.
\end{remark}

\vskip .5cm

\noindent \textbf{The updates for $(a,b,\alpha,\beta,\phi)$.} Note that the updates for $a,b,\alpha,\beta,\phi$ are decoupled. Indeed, \eqref{eq:algorithm} yields
\begin{equation*}
	\begin{cases}
	a^{n}(t)+\tau \left(\int_\Omega \rho^n(x,t) \zeta_i(x)dx\right)_i \in a^{n+1}(t)+\tau \bK^{-1} \partial_a U_1^*(a^{n+1}(t))\\
	b^{n}+\tau \left(\int_\Omega \rho^n(x,1) \zeta_i(x)dx\right)_i \in b^{n+1}+ \tau \bS^{-1} \partial_b V_1^*(b^{n+1})\\
	\alpha^n(x,t)+\tau \rho^n(x,t) \in \alpha^{n+1}(x,t)+\tau \partial_\alpha U_0^*(\alpha^{n+1}(x,t))\\
	\beta^n(x)+\tau \rho^n(x,1) \in \beta^{n+1}(x)+\tau \partial_\beta V^*_0(\beta^{n+1}(x))\\
	\phi^n(x,t)+\tau \mathcal{L}^{-1}(\rho^n_t+\nabla \cdot m^n,0,0)=\phi^{n+1}(x,t)
	\end{cases}
\end{equation*}
To update $a$, we need to solve an $r\times r$ system for every fixed $t$. Next, to update $b$ we need to solve an $r\times r$ system. Once $r$ is fixed the sizes of these systems do not depend on the mesh.

Next, we observe that
\begin{equation*}
\begin{split}
U_0^*(\alpha)=&\int_0^1 \int_\Omega F_0^*(t,x,\alpha(x,t)) dxdt,\quad \partial_\alpha U_0^*(\alpha(x,t))=\partial_\alpha F_0^*(t,x,\alpha(x,t)),\\
V_0^*(\beta)=&\int_\Omega G_0^*(x,\beta(x)) dx,\quad \partial_\beta V_0^*(\beta(x))=\partial_\beta G_0^*(x,\beta(x)),
\end{split}
\end{equation*}
where
\begin{equation*}
F_0^*(t,x,\alpha)=\sup_\rho \alpha \rho -F_0(t,x,\rho),\quad G_0^*(x,\beta)=\sup_\rho \beta \rho -G_0(x,\rho).
\end{equation*}
Therefore, the updates for $\alpha,\beta$ correspond to decoupled one-dimensional proximal steps; that is,
\begin{equation*}
\begin{cases}
\alpha^{n+1}(x,t) \in \arg\min_{\alpha} F_0^*(t,x,\alpha) +\frac{|\alpha-\alpha^n(x,t)-\tau \rho^n(x,t)|^2}{2\tau}\\
\beta^{n+1} \in \arg\min_\beta G_0^*(x,\beta) +\frac{|\beta-\beta^n(x)-\tau \rho^n(x,1)|^2}{2\tau}
\end{cases}
\end{equation*}
Therefore, the updates for $\alpha,\beta$ can be efficiently performed in parallel yielding linear-in-grid computational cost. Finally, recalling the definition of $\mathcal{L}$, we obtain that to update $\phi$ we need to solve a space-time elliptic equation
\begin{equation*}
	\begin{cases}
	\Delta_{t,x} \phi=\Delta_{t,x} \phi^n+\tau (\rho^n_t+\nabla \cdot m^n)\\
	\phi(x,0)-\phi_t(x,0)=\phi^n(x,0)-\phi^n_t(x,0)\\
	\phi(x,1)+\phi_t(x,1)=\phi^n(x,1)+\phi^n_t(x,1)
	\end{cases}
\end{equation*}
This step can be efficiently performed via Fast Fourier Transform (FFT).

\vskip .5cm

\noindent\textbf{The updates for $(\rho,m)$.} The resolvent operator $J_{\sigma \partial \Xi}$ is the proximal operator $\operatorname{prox}_{\sigma \Xi}$. Therefore, $(\rho,m)$ updates reduce to an optimization problem
\begin{equation*}
	\begin{split}
	\inf_{\rho,m} &\int_0^1 \int_\Omega \rho L \left(x,\frac{m}{\rho}\right)dxdt+\mathbf{1}_{\rho(x,0)=\rho_0(x)}+\mathbf{1}_{\rho\geq 0}+\mathbf{1}_{m<<\rho}\\
	&+\frac{1}{2\sigma}\int_{\Omega}\left(\rho(x,0)-\rho^n(x,0)+\sigma\tilde{\phi}^{n+1}(x,0)\right)^2dx\\
	&+\frac{1}{2\sigma}\int_{\Omega}\left(\rho(x,1)-\rho^n(x,1)+\sigma \sum_i \tilde{b}_i^{n+1} \zeta_i(x)+\sigma \tilde{\beta}^{n+1}(x)-\sigma \tilde{\phi}^{n+1}(x,1) \right)^2dx\\
	&+\frac{1}{2\sigma} \int_0^1 \int_\Omega \left(\rho(x,t)-\rho^n(x,t) +\sigma \sum_i \tilde{a}^{n+1}(t)\zeta_i(x)+\sigma \tilde{\alpha}^{n+1}(x,t)+\sigma \tilde{\phi}^{n+1}_t(x,t) \right)^2dxdt\\
	&+\frac{1}{2\sigma} \int_0^1 \int_\Omega \left| m(x,t)-m^n(x,t)+\sigma \nabla \tilde{\phi}^{n+1}(x,t)\right|^2 dxdt
	\end{split}
\end{equation*}
Again, we obtain decoupled one-dimensional optimization problems
\begin{equation*}
\begin{cases}
(\rho^{n+1}(x,t),m^{n+1}(x,t)) \in \arg\min\limits_{\rho,m} \rho L \left(x,\frac{m}{\rho}\right)+\mathbf{1}_{\rho\geq 0}+\mathbf{1}_{m<<\rho} \\
+\frac{\left|\rho-\rho^n(x,t) +\sigma \sum_i \tilde{a}^{n+1}(t)\zeta_i(x)+\sigma \tilde{\alpha}^{n+1}(x,t)+\sigma \tilde{\phi}^{n+1}_t(x,t)\right|^2}{2\sigma }\\
+\frac{\left| m-m^n(x,t)+\sigma \nabla \tilde{\phi}^{n+1}(x,t) \right|}{2\sigma}\\
\rho^{n+1}(x,0) \in \arg\min\limits_{\rho} \mathbf{1}_{\rho=\rho_0(x)}+\frac{\left|\rho-\rho^n(x,0)+\sigma\tilde{\phi}^{n+1}(x,0)\right|^2}{2\sigma}\\
\rho^{n+1}(x,1) \in \arg\min\limits_{\rho} \frac{\left|\rho-\rho^n(x,1)+\sigma \sum_i \tilde{b}_i^{n+1} \zeta_i(x)+\sigma \tilde{\beta}^{n+1}(x)-\sigma \tilde{\phi}^{n+1}(x,1)\right|^2}{2\sigma}\\
\end{cases}
\end{equation*}

\section{A class of non-potential MFG with density constraints}\label{sec:kernels}

Here discuss an instance of \eqref{eq:main} that is non-potential and incorporates pointwise density constraints for the agents. We illustrate that our method handles mixed couplings in an efficient manner. Assume that
\begin{equation}\label{eq:fg_lin_nonlocal}
\begin{split}
f_1(t,x,z)=&z,\quad f_0(t,x,z)=\partial_z \mathbf{1}_{\underline{h}(x,t)\leq z \leq \bar{h}(x,t)}\\
g_1(x,z)=&z,\quad g_0(x,z)=\partial_z \mathbf{1}_{\underline{e}(x)\leq z \leq \bar{e}(x)}+g(x)
\end{split}
\end{equation}
Functions $\underline{h},~\bar{h}\geq 0$ and $\underline{e},\bar{e}\geq 0$ are density constraints; that is, the solution to the MFG problem must satisfy the hard constraints $\underline{h}(x,t)\leq \rho(x,t) \leq \bar{h}(x,t)$, $(x,t)\in \Omega \times (0,1)$ and $\underline{e}(x)\leq \rho(x,1)\leq \bar{e}(x)$, $x\in \Omega$. Next, $g$ is a terminal cost function.
\begin{remark}
We can model static and dynamic obstacles in this framework. Indeed, assume that $\Omega_t \subset \Omega$ is a dynamic obstacle and set $\underline{h}(x,t)=\bar{h}(x,t)=\chi_{\Omega_t}(x)$. Then the hard constraint $\underline{h}(x,t)\leq \rho(x,t) \leq \bar{h}(x,t)$ is equivalent to $\mathrm{supp}\rho(\cdot,t) \cap \Omega_t=\emptyset$, which means that there are no agents in $\Omega_t$. One can also use the lower bounds on $\rho$ to maintain a minimal fraction of agents at specific locations.
\end{remark}

From \eqref{eq:U}, we obtain
\begin{equation*}
	\begin{split}
	U_1(c)=&\int_0^1 \int_\Omega \frac{\left(\sum_{j=1}^r c_j(t) \zeta_j(x)\right)^2}{2}dxdt=\frac{1}{2} \int_0^1 \sum_{j=1}^r c_j^2(t)dt\\
	V_1(w)=&\int_\Omega \frac{\left(\sum_{j=1}^r w_j \zeta_j(x)\right)^2}{2}dx=\frac{1}{2} \sum_{j=1}^r w_j^2\\
	U_{0}(\rho)=&\int_0^1 \int_{\Omega} \mathbf{1}_{\underline{h}(x,t)\leq \rho(x,t) \leq \bar{h}(x,t)} dxdt\\
	V_{0}(\mu)=& \int_{\Omega} \mathbf{1}_{\underline{e}(x)\leq \mu(x) \leq \bar{e}(x)}+g(x)\mu(x) dx
	\end{split}
\end{equation*}
Furthermore, the dual functions are
\begin{equation*}
\begin{split}
U_1^*(a)=&\frac{1}{2} \int_0^1 \sum_{j=1}^r a_j^2(t)dt,\quad V_1^*(b)=\frac{1}{2} \sum_{j=1}^r b_j^2\\
U_{0}^*(\alpha)=&\int_0^1 \int_{\Omega} \max \left\{\alpha(x,t) \underline{h}(x,t), \alpha(x,t) \bar{h}(x,t) \right\} dxdt\\
V_{0}^*(\beta)=& \int_{\Omega} \max \left\{\left(\beta(x)-g(x)\right) \underline{e}(x), \left(\beta(x)-g(x)\right) \bar{e}(x) \right\} dx\\
\end{split}
\end{equation*}
Accordingly, the algorithm \eqref{eq:algorithm} reduces to
\begin{equation*}
\begin{cases}
a^{n+1}=\left(\mathbf{I}+\tau \bK^{-1}\right)^{-1}\left(a^n+\tau \left(\rho^n(x,t)\zeta_j(x)\right)_j\right)\\
\alpha^{n+1}(x,t)=\min \left\{ \max \left\{ 0,\alpha^n(x,t)+\tau \rho^n(x,t)-\tau \bar{h}(x,t) \right\} , \alpha^n(x,t)+\tau \rho^n(x,t)-\tau \underline{h}(x,t) \right\}\\
\beta^{n+1}(x)=\min \left\{ \max \left\{g(x),\beta^n(x)+\tau \rho^{n}(x,1)-\tau \bar{e}(x)\right\}, \beta^n(x)+\tau \rho^{n}(x,1)-\tau \underline{e}(x) \right\}\\
\Delta_{t,x} \phi^{n+1}=\Delta_{t,x} \phi^n+\tau (\rho^n_t+\nabla \cdot m^n)\\
\phi^{n+1}(x,0)-\phi^{n+1}_t(x,0)=\phi^n(x,0)-\phi^n_t(x,0)\\
\phi^{n+1}(x,1)+\phi^{n+1}_t(x,1)=\phi^n(x,1)+\phi^n_t(x,1)\\
\left(\tilde{a}^{n+1},\tilde{\alpha}^{n+1},\tilde{\beta}^{n+1},\tilde{\phi}^{n+1}\right)=2\left(a^{n+1},\alpha^{n+1},\beta^{n+1},\phi^{n+1}\right)-\left(a^{n},\alpha^{n},\beta^{n},\phi^n\right)\\
(\rho^{n+1}(x,t),m^{n+1}(x,t)) \in \arg\min\limits_{\rho,m} \rho L \left(x,\frac{m}{\rho}\right) +\mathbf{1}_{\rho\geq 0}+\mathbf{1}_{m<<\rho}\\
+\frac{\left|\rho-\rho^n(x,t) +\sigma \sum_i \tilde{a}^{n+1}(t)\zeta_i(x)+\sigma \tilde{\alpha}_1^{n+1}(x,t)+\sigma \tilde{\alpha}_2^{n+1}(x,t)+\sigma \tilde{\phi}^{n+1}_t(x,t)\right|^2}{2\sigma }\\
+\frac{\left| m-m^n(x,t)+\sigma \nabla \tilde{\phi}^{n+1}(x,t) \right|}{2\sigma}\\
\rho^{n+1}(x,0)=\rho_0(x)\\
\rho^{n+1}(x,1)=\rho^n(x,1)-\sigma \tilde{\beta}^{n+1}(x)+\sigma \tilde{\phi}^{n+1}(x,1)
\end{cases}
\end{equation*}

We have that $\bK=(k_{pq})$ where $k_{pq}=\int_{\Omega^2} K(x,y) \zeta_{p}(x)\zeta_q(y)dxdy$. Therefore, $\bK$ may not be symmetric if $K$ is not. In this case, \eqref{eq:main} is non-potential. Nevertheless, if $K$ is monotone then such is $\bK$, and our methods apply. Below we discuss a class of non-symmetric interactions that are monotone but non-symmetric. For $\delta_-,\delta_+>0$ consider
\begin{equation*}
	\gamma_{\delta_-,\delta_+}(x)=e^{-\frac{x^2}{2\delta_-^2}} \chi_{x<0}+e^{-\frac{x^2}{2\delta_+^2}} \chi_{x\geq 0},\quad x\in \mathbb{R}
\end{equation*}
The cosine transform of $\gamma$ is
\begin{equation*}
	\int_{\mathbb{R}} \cos \left(2 i \pi \zeta x\right) \gamma(x)dx=\sqrt{\frac{\pi}{2}} \left(\delta_- e^{-2\pi^2 \zeta^2 \delta_-^2} +\delta_+ e^{-2\pi^2 \zeta^2 \delta_+^2}\right)>0,\quad \zeta \in \mathbb{R}.
\end{equation*}
Therefore, $K(x,y)=\gamma(y-x)$ is a monotone kernel. Therefore, for $\delta_-,\delta_+ \in \mathbb{R}_+^d$
\begin{equation*}
	K_{\delta_-,\delta_+}(x,y)=\prod_{i=1}^d \gamma_{\delta_{i,-},\delta_{i,+}}(y_i-x_i)
\end{equation*}
is a monotone kernel. Furthermore, for any non-singular linear transformation $Q$ we have that $K_{\delta_-,\delta_+}(Q^{-1}x,Q^{-1}y)$ is a monotone kernel. Therefore, for a basis $\nu=\{\nu_1,\nu_2,\cdots,\nu_d\} \subset \mathbb{R}^d$ we have that
\begin{equation}\label{eq:nonsym_Gaussian}
	K_{\nu,\delta_-,\delta_+}(x,y)=K_{\delta_-,\delta_+}(Q^{-1}x,Q^{-1}y)=\prod_{i=1}^d \gamma_{\delta_{i,-},\delta_{i,+}}(y'_i-x'_i)
\end{equation} 
is a monotone kernel, where $Q=(\nu_1 ~\nu_2~ \cdots ~\nu_d)$ is the coordinates transformation matrix and $x'=Q^{-1}x$ are the coordinates in $\nu$.

Kernels in \eqref{eq:nonsym_Gaussian} model interactions that have different strengths of repulsion along lines parallel to $\nu_i$. Moreover, these interactions are not symmetric as they depend on the sign of $y_i'-x_i'$ that tells us whether $y$ is in the front or back of $x$ relative to $\nu_i$. We can think of crowd motion models where people mostly pay attention to the crowd in front of them.

\section{Numerical experiments}\label{sec:numerics}

In this section, we present three sets of numerical examples for MFG with mixed couplings using Algorithm \eqref{eq:algorithm}. We take $\Omega \times [0,T] = [-1,1]^2 \times [0,1]$, with a uniform space-grid $N_x= 64$ and a uniform time-grid $N_t = 32$ in all examples. 
  \subsection{Density Splitting with Asymmetric Kernel}
  We consider a MFG problem where the density splits into $8$ parts at final time.
  For this non-potential MFG with density constraints \eqref{eq:fg_linear_nonlocal}, we set
  \begin{equation*}
  \begin{split}
  \rho_0(x) &= \mathcal{N}([0,0],0.1) \\
  g(x) & =\frac{1}{2} \sum_{j=1}^8 \left(1 - \exp\left(20\|x-x_j\|^2\right)\right), \quad\text{where}\; x_j = 0.75 \left[\sin{\frac{2\pi j}{8}},\cos{\frac{2\pi j}{8}} \right]\\
  \bar{e}(x) & = 4,
  \end{split}
  \end{equation*}
  where  $\mathcal{N}(x,\sigma_G)$ is the density of a homogeneous normal distribution centered at $x$ with variance $\sigma_G$.
  As for the Gaussian type Kernel $K(x,y)$ in \eqref{eq:nonsym_Gaussian}, we choose the three following set-up:

  \begin{itemize}
  	\item Case A, symmetric kernel
  	\begin{equation*}
  	\begin{split}	
  	&K(x,y) = 4\exp\left(-\frac{\|x-y\|^2 }{2\delta^2} \right)\\
  	&\delta = \delta_{i,-}=\delta_{i,+}=0.1 \;\text{for} \; i=1,2
  	\end{split}
  	\end{equation*}
  	\item Case B, asymmetric kernel
  	\begin{equation*}
  	\begin{split}	
  	& K(x,y) = 4 \gamma_{\delta_{-},\delta}(x_1-y_1)\gamma_{\delta,\delta}(x_2-y_2)\\
  	&\delta = 0.1, \delta_{-}=0.4.
  	\end{split}
  	\end{equation*}	
  	\item Case C, asymmetric kernel (coordination transform) 
  	\begin{equation*}
  	\begin{split}
  	&K(x,y) = 4\exp\left(-\frac{(x-y)^T Q(x-y) }{2\delta^2} \right)\\
  	&	Q = \begin{pmatrix}
  	\frac{1}{\delta^2} & 	\frac{c}{\delta^2} \\
  	\frac{c}{\delta^2} & 	\frac{1}{\delta^2}
  	\end{pmatrix}, c=0.95, \delta =0.1
  	\end{split}
  	\end{equation*}	
  \end{itemize}
  \begin{figure} [!htbp]
  	\centering
  	\begin{subfigure}{\linewidth}
  		\includegraphics[width=1.0\textwidth,trim=0 0 0 0, clip=false]{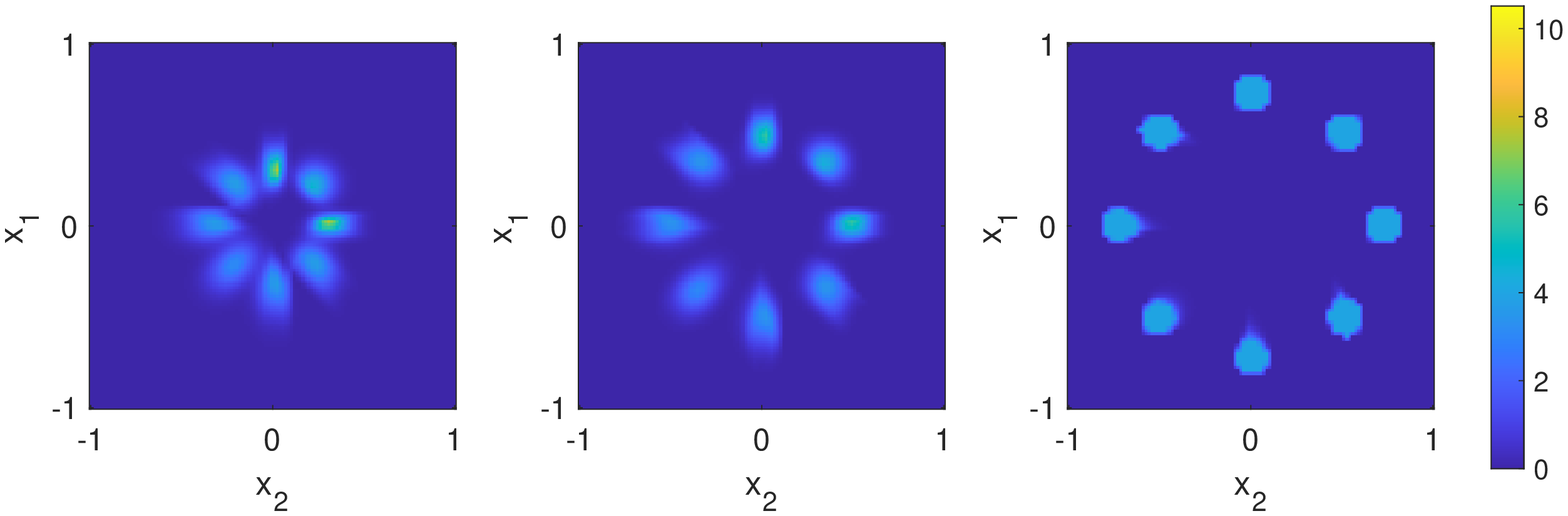}
  		\vspace{-0.7cm}
  		\caption{Symmetric kernel}
  	\end{subfigure}
  	\begin{subfigure}{\linewidth}
  		\includegraphics[width=1.0\textwidth,trim=0 0 0 0, clip=false]{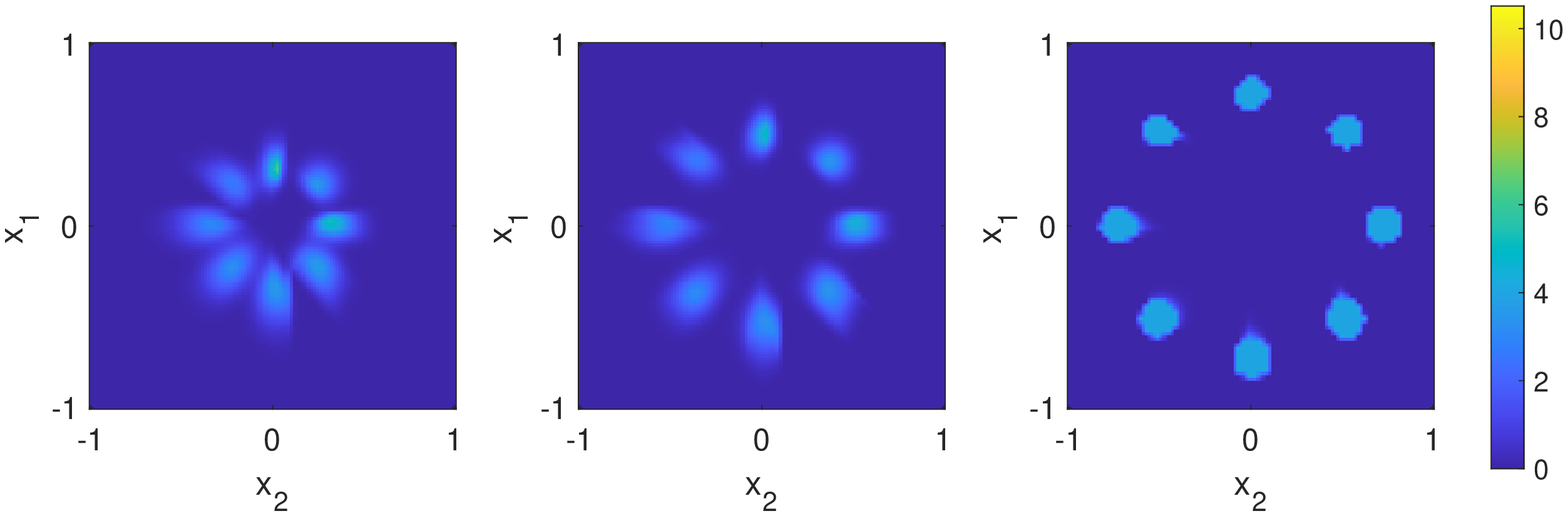}
  		\vspace{-0.7cm}
  		\caption{Asymmetric kernel}
  	\end{subfigure}
  	\begin{subfigure}{\linewidth}
  		\includegraphics[width=1.0\textwidth,trim=0 0 0 0, clip=false]{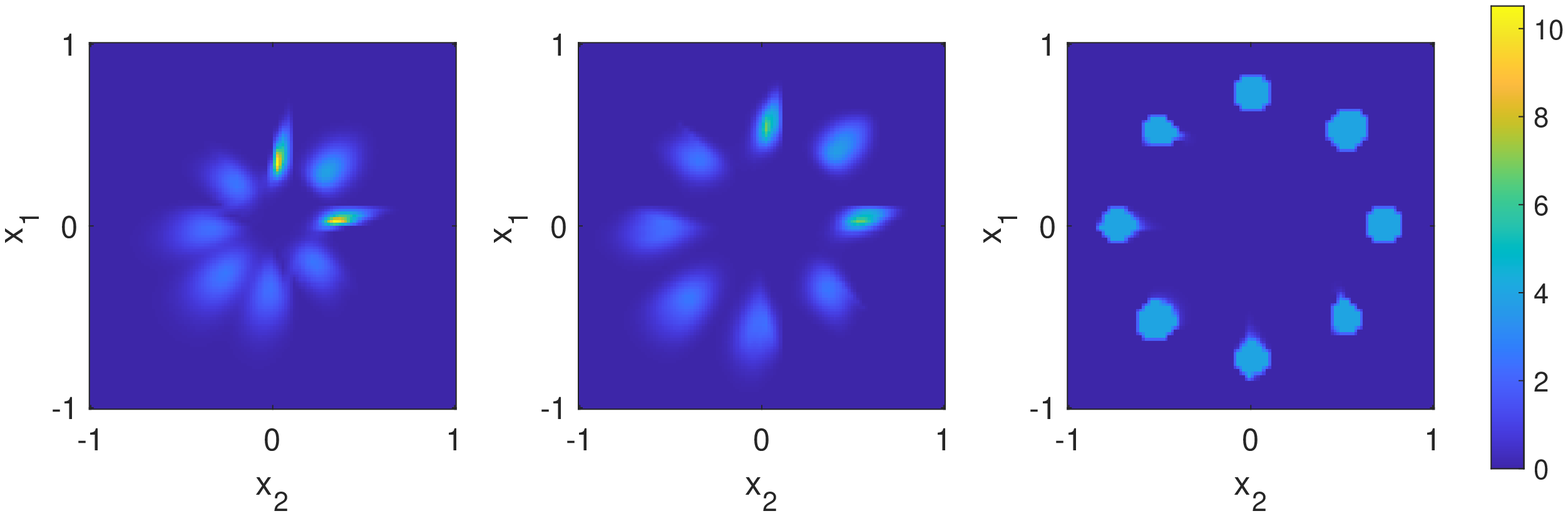}
  		\vspace{-0.7cm}
  		\caption{Asymmetric kernel (coordination transform) }
  	\end{subfigure}
  	\caption{MFG solution $\rho(x,0.3), \rho(x,0.6), \rho(x,1)$ for density splitting examples.}
  	\label{fig:example1}
  \end{figure} 
  The results are shown in Figure \ref{fig:example1}. As we can see nonlocal kernel affect how the density moves and have different final distributions.
  In case A, a symmetric kernel leads to an even splitting of the initial density. 
  Comparing case A and B, we see that the large $\delta_{1,-}$ causes density to favor a motion towards $x_{1,-}$ direction. As a result, the final density has more concentration in $x_{1,-}$ domain.  As for a comparison between A and C, we see that agents in case C have a preference to move along $x_1 = x_2$ direction, which is consistent with the shape of the kernel $K(x,0)$.

  \subsection{Static Obstacles Modeled with Density constraint}
  Here we provide a MFG problem where the the density moves while avoiding the obstacles, which is modeled using density constraint. We also include a small local interaction term in this example:
  \begin{equation*}
  \begin{split}
  & f_0(t,x,z)=\partial_z \mathbf{1}_{\underline{h}(x,t)\leq z \leq \bar{h}(x,t)}	+ \epsilon \log{z}, \quad \epsilon = 0.01, \\
  &\bar{h}(x,t) = 0, \bar{e}(x) = 0 \;\text{for}\; x\in \Omega_{obs},\\
  & K(x,y) = 4 \gamma_{\delta_{-},\delta}(x_1-y_1)\gamma_{\delta,\delta}(x_2-y_2), \quad \delta = 0.1, \delta_{-}=0.4
  \end{split}
  \end{equation*}	
  As for the initial-terminal conditions, we have 
  \begin{equation*}
  \begin{split}
  & \rho_0(x) = \frac{1}{2}\mathcal{N}([-0.8,0.5],0.1) + \frac{1}{2}\mathcal{N}([-0.8,0.5],0.1)\\
  & g(x) = x_1^2 + \left(x_2-0.85\right)^2 -2e^{-10\left(x_2-0.75\right)^2}.
  \end{split}
  \end{equation*}	
  Under above setup, we consider the following $2$ cases:
  case A, $\bar{h}(x,t) = 20$; case B, $\bar{h}(x,t) = 10$ for all $(x,t), x\notin \Omega_{obs}$.
  \begin{figure} [!htbp]
  	\centering
  	\begin{subfigure}{\linewidth}
  		\includegraphics[width=1.0\textwidth,trim=0 0 0 0, clip=false]{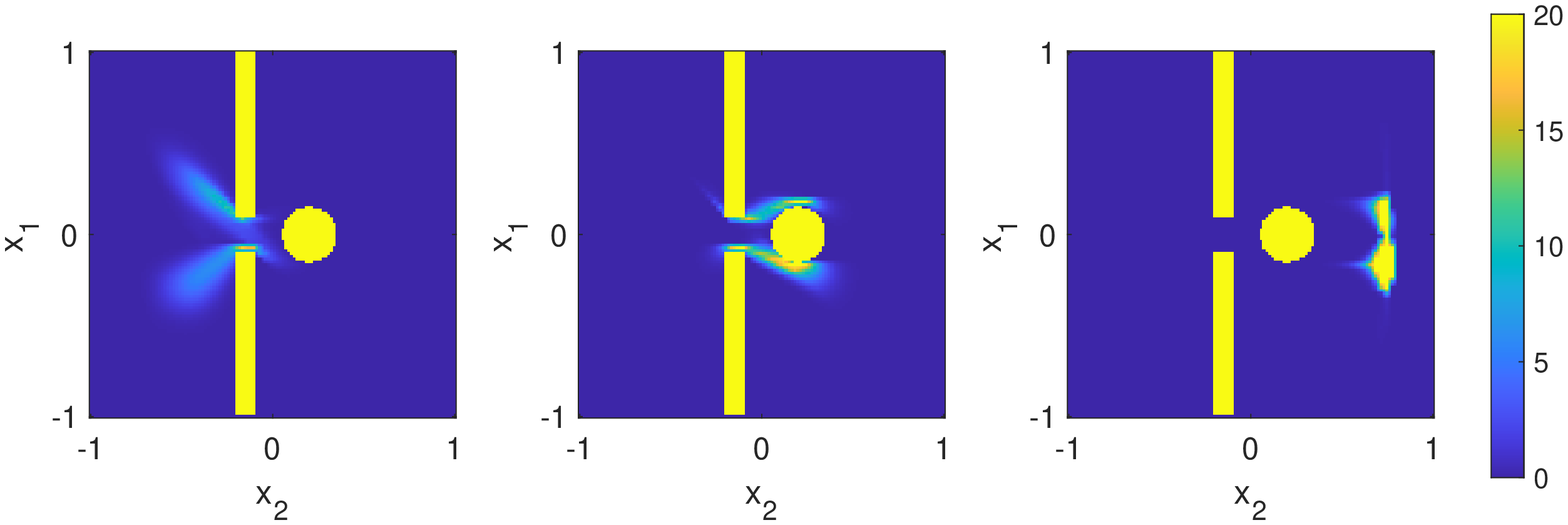}
  		\vspace{-0.7cm}
  		\caption{$\bar{h}(x,t) = 20$}
  	\end{subfigure}
  	\begin{subfigure}{\linewidth}
  		\includegraphics[width=1.0\textwidth,trim=0 0 0 0, clip=false]{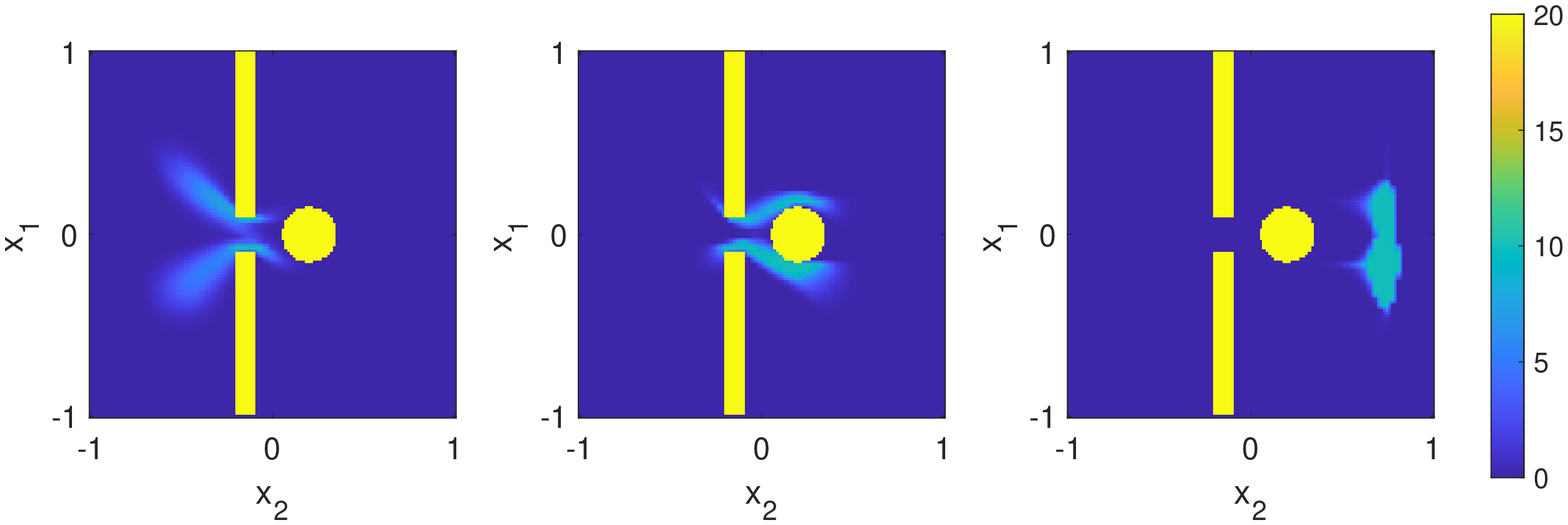}
  		\vspace{-0.7cm}
  		\caption{ $\bar{h}(x,t) = 10$}
  	\end{subfigure}
  	
  	\caption{MFG solution $\rho(x,0.3), \rho(x,0.6), \rho(x,1)$ for static obstacles examples, where the obstacle (yellow) is located at $\Omega_{obs} = \left\{\|x-[0,0.2]\|^2\leq 0.15^2\right\} \cup \left\{ |x_1|\geq 0.1, |x_2+0.15|\leq 0.05 \right\}  $ }
  	\label{fig:example2}
  \end{figure} 
  The numerical results are shown in Figure \ref{fig:example2}. As we see the density moves from left to the right and avoids both the rectangle and round obstacles. It avoids the round obstacles via an uneven splitting, which is caused by the asymmetric kernel. Comparing Case A and B, we see that the density constraint make the agents spread more.
  
  \subsection{Dynamic Obstacles Modeled with Density constraint}
\begin{figure} [!htbp]
  	\centering
  	\begin{subfigure}{\linewidth}
  		\includegraphics[width=1.0\textwidth,trim=0 0 0 0, clip=false]{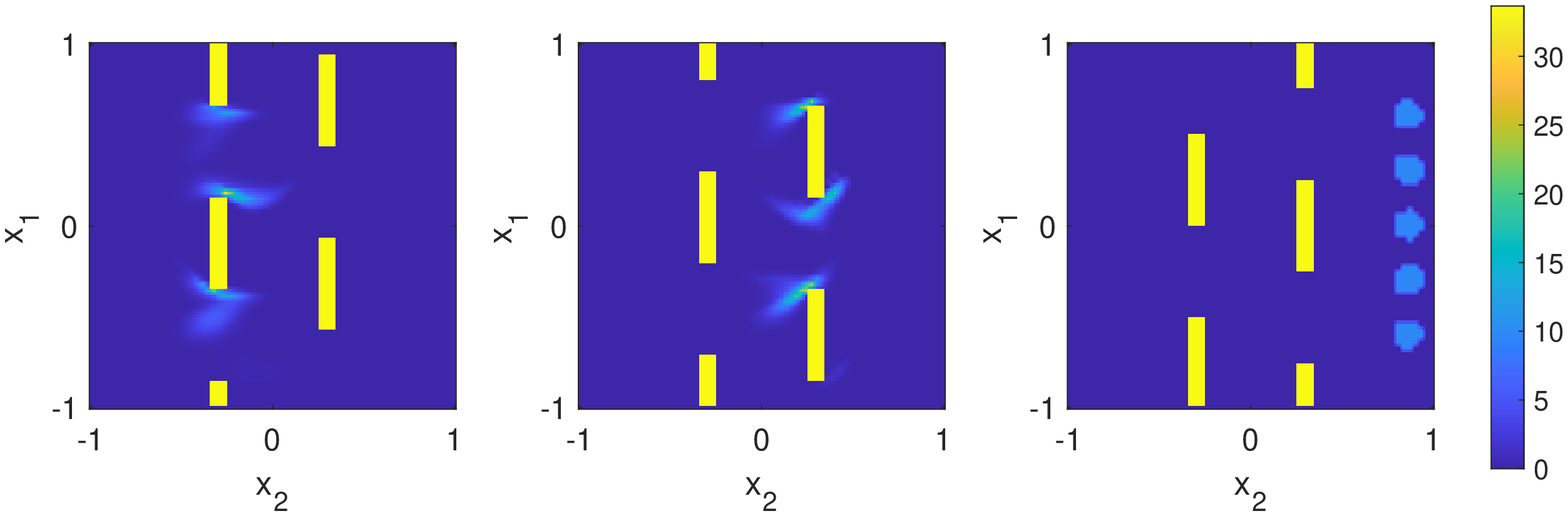}
  		\vspace{-0.7cm}
  		\caption{$\underline{e}(x)= c_e$}
  	\end{subfigure}
  	\begin{subfigure}{\linewidth}
  		\includegraphics[width=1.0\textwidth,trim=0 0 0 0, clip=false]{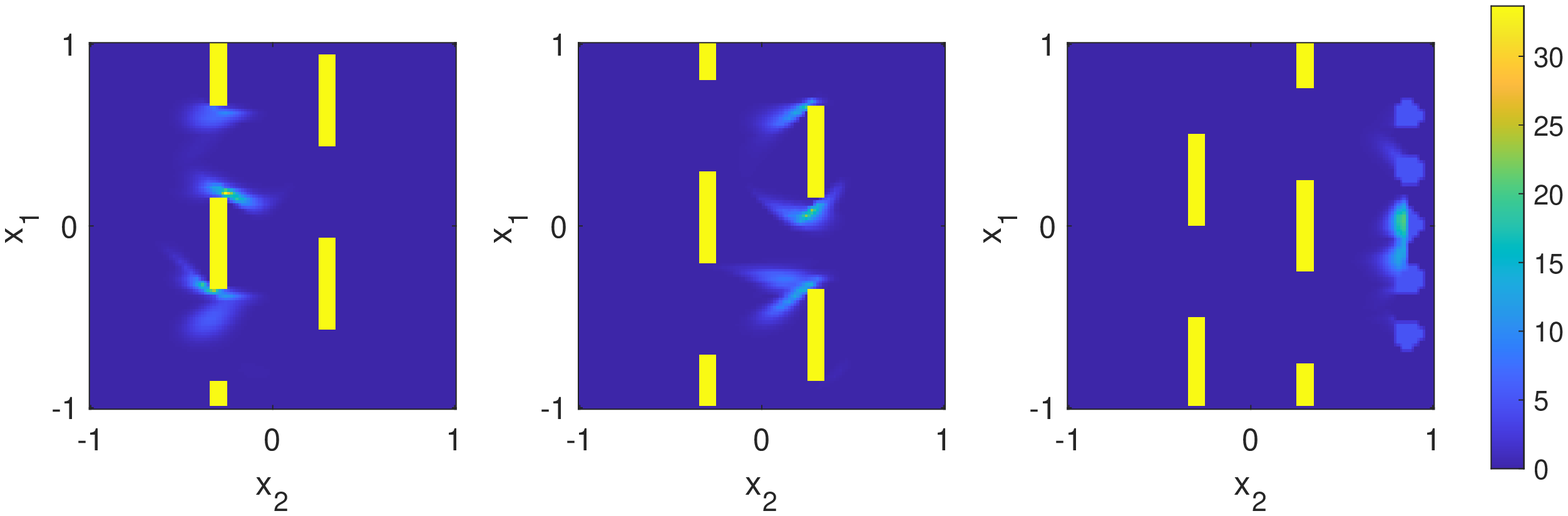}
  		\vspace{-0.7cm}
  		\caption{ $\underline{e}(x)= \frac{1}{2}c_e$}
  	\end{subfigure}  	
  	\caption{MFG solution $\rho(x,0.3), \rho(x,0.6), \rho(x,1)$ for dynamic obstacles examples.}
  	\label{fig:example3}
  \end{figure}

  In the last example, we model an Optimal-Transport-like problem with dynamic obstacles via density constraint. Our mean field game system \label{eq:fg_linear_nonlocal} is as follows:
  \begin{equation*}
  \begin{split}
  &	K(x,y) = \exp\left(-\frac{\|x-y\|^2 }{2\delta^2} \right),\quad 	\delta=0.1 \\
  &	\rho_0(x) = \sum_{j=1}^5 \frac{1}{5}\mathcal{N}(x_j,0.05), \text{for}\; x_j = [-0.9+0.3j,-0.85]\\
  & g(x)= x_1^2 + 5|x_2-0.85|^{1.5}.
  \end{split}
  \end{equation*}
  We use the density constraint $\bar{h}(x,t)$ to model 4 rectangles moving vertically. As for the density constraint at final time, we choose $\underline{e}(x)$ to be exactly a density distribution. Specifically,
  \begin{equation*}
  \underline{e}(x)= c_e,\; \text{for}\; \|x-x_j\|\leq 0.08,\; \text{for} \; x_j= [-0.9+0.3j,0.85], j = 1...5,
  \end{equation*}
  where $c_e$ is a constant that normalized $\underline{e}(x)$. This setup is equivalent to specifying the final density distribution $\rho(x,1)=\underline{e}(x)$. This is case A that is shown in Figure \ref{fig:example3}.
  We also relax the density constraint at $t=1$, by setting 
  \begin{equation*}
  \underline{e}(x)=\frac{1}{2} c_e,\; \text{for}\; \|x-x_j\|\leq 0.08,\; \text{for} \; x_j= [-0.9+0.3j,0.85], j = 1...5.
  \end{equation*}
  That it, we decrease the lower bound of $\rho(x,1)$ by half. This is Case B shown in Figure \ref{fig:example3}. As we can see, unlike case A, the agents do not completely fill the support of $\underline{e}(x)$. The aspect of modeling dynamic obstacles also works well as agents avoid the prohibited regions.

\section*{Appendix}

\subsection*{Proof of Proposition \ref{prp:Psi}}
Calculating the first variation with respect to $\rho,m$ in \eqref{eq:Psi}, we obtain
\begin{equation*}
	\begin{cases}
	L \left(t,x,\frac{m}{\rho}\right)-\frac{m}{\rho} \nabla_v L \left(t,x,\frac{m}{\rho}\right)+\alpha(x,t)+\sum_{i=1}^r a_i(t) \zeta_i(x)+\phi_t=0\\
	\nabla_v L \left(t,x,\frac{m}{\rho}\right)+\nabla \phi=0\\
	\beta(x)+\sum_{i=1}^r b_i \zeta_i(x)-\phi(x,1)=0\\
	\end{cases}
\end{equation*}
Additionally, $\rho(x,0)=\rho_0(x)$ as $\Xi(\rho,m)<\infty$. From the properties of the Legendre transform, we obtain that
\begin{equation*}
	\begin{cases}
	\frac{m}{\rho}=-\nabla_p H(t,x,\nabla \phi)\\
	H(t,x,\nabla \phi)= \frac{m}{\rho} \nabla_v L \left(t,x,\frac{m}{\rho}\right)-L \left(t,x,\frac{m}{\rho}\right)
	\end{cases}
\end{equation*}
which yields \eqref{eq:main_r_alpha}. Rest of the proof readily follows from the convex duality relation between $\Xi$ and $\Xi^*$.

\subsection*{Proof of Proposition \ref{prp:M}}
The components of $C^*$ corresponding to variables $a,b,\alpha,\beta$ are straightforward to calculate as they are in $L^2$ spaces. As for the component corresponding to $\phi$, we have to find $h=h(\rho,m)$ such that for all $\phi$ one has that
\begin{equation*}
\begin{split}
\langle h, \phi \rangle_{H^1}=& -\int_0^1 \int_\Omega \phi_t \rho +\nabla \phi \cdot mdxdt+\int_{\Omega} \phi(x,1) \rho(x,1)dx-\int_{\Omega} \phi(x,0) \rho(x,0)dx\\
=&\int_0^1 \int_\Omega \phi \left(\rho_t +\nabla \cdot m\right) dxdt
\end{split}
\end{equation*}
We have that
\begin{equation*}
	\begin{split}
	\la h,\phi \ra_{H^1}=&\int_\Omega \phi(x,0) h(x,0) dx+\int_\Omega \phi(x,1) h(x,1) dx+\int_0^1 \int_\Omega \nabla_{t,x} \phi \cdot \nabla_{t,x} h dx dt\\
	=&-\int_0^1 \int_\Omega \phi \Delta_{t,x} h dxdt+\int_\Omega \phi(x,1) \left(h(x,1)+h_t(x,1)\right)dx\\
	&+\int_\Omega \phi(x,0) \left(h(x,0)-h_t(x,0)\right) dx
	\end{split}
\end{equation*}
Therefore $h$ must satisfy the conditions
\begin{equation}\label{eq:h}
\begin{cases}
\Delta_{t,x} h=-\left(\rho_t + \nabla \cdot m\right)\\
h(x,0)-h_t(x,0)=0\\
h(x,1)+h_t(x,1)=0
\end{cases}
\end{equation}
and we obtain \eqref{eq:C*}. Next, we prove \eqref{eq:C*rho_alpha}. We have that
\begin{equation*}
	U_0^*(\alpha)=\int_0^1 \int_\Omega F_0^*(t,x,\alpha(x,t)) dxdt,\quad \partial_\alpha U_0^*(\alpha(x,t))=\partial_\alpha F_0^*(t,x,\alpha(x,t)),
\end{equation*}
where $F_0^*(t,x,\alpha)=\sup_\rho \alpha \rho -F_0(t,x,\rho)$. Therefore, the $\alpha$-entry inclusion in \eqref{eq:C*rho_alpha} is equivalent to
\begin{equation*}
	\rho(x,t) \in \partial_\alpha F_0^*(t,x,\alpha(x,t)) \Longleftrightarrow \alpha(x,t) \in \partial_\rho F_0(t,x,\rho)=f_0(t,x,\rho)
\end{equation*}
Similarly, the $\beta$-entry inclusion is equivalent to
\begin{equation*}
\rho(x,1) \in \partial_\alpha G_0^*(x,\beta(x)) \Longleftrightarrow \beta(x) \in \partial_\rho G_0(x,\rho(x,1))=g_0(x,\rho(x,1))
\end{equation*}
Next, the $\phi$-entry inclusion means that $h=0$ in \eqref{eq:h} that is equivalent to
\begin{equation*}
	\rho_t+\nabla \cdot m=0
\end{equation*}
The $a$-entry inclusion in \eqref{eq:C*rho_alpha} is
\begin{equation*}
	\left(\int_\Omega\rho(x,t)\zeta_i(x) dx\right)_{i=1}^r \in \bK^{-1} \partial_a U_1^*(a) \Longleftrightarrow  \bK \left(\int_\Omega\rho(x,t)\zeta_i(x)dx\right)_{i=1}^r \in \partial_a U_1^*(a)
\end{equation*}
Applying the properties of the Legendre transform again, we obtain that this previous inclusion is equivalent to
\begin{equation}\label{eq:a_entry}
	a \in \partial_c U_1\left(\bK \left(\int_\Omega \rho(x,t)\zeta_i(x)\right)_{i=1}^r \right)
\end{equation}
On the other hand, we have that
\begin{equation*}
\begin{split}
\partial_{c_i} U_1(c)=\int_0^1 \int_\Omega f_1\left(t,x,\sum_{p=1}^r c_p(t) \zeta_p(x)\right) \zeta_i(x)dxdt
\end{split}
\end{equation*}
and therefore
\begin{equation*}
\begin{split}
&\partial_{c_i} U_1\left(\bK \left(\rho(x,t)\zeta_i(x)\right)_{i=1}^r \right)\\
=&\int_0^1 \int_\Omega f_1\left(t,x,\sum_{p=1}^r \sum_{q=1}^r k_{pq} \int_{\Omega} \rho(y,t) \zeta_q(y)dy \zeta_p(x)\right) \zeta_i(x)dxdt\\
=&\int_0^1 \int_\Omega f_1\left(t,x, \int_{\Omega} K_r(x,y)\rho(y,t) dy \right) \zeta_i(x)dxdt
\end{split}
\end{equation*}
Hence, \eqref{eq:a_entry} is equivalent to the first equation in \eqref{eq:consistency}. The derivation for the $b$-entry in \eqref{eq:C*rho_alpha} is similar.

\subsection*{Proof of Theorem \ref{thm:inclusion}}
The equivalence of \eqref{eq:inclusion_ab_phi} and \eqref{eq:main_r_alpha}-\eqref{eq:consistency} is simply a combination of assertions in Propositions \ref{prp:Psi} and \ref{prp:M}. Furthermore, assume that $c \mapsto \partial_c U_1(\bK c)$ and $w\mapsto \partial_w V_1(\bS w)$ are maximally monotone. We have that
\begin{equation*}
	\left(\partial_c U_1\left(\bK c\right)\right)^{-1}=\bK^{-1} \partial_a U_1^*(a),\quad \left(\partial_w V_1\left(\bS w\right)\right)^{-1}=\bS^{-1} \partial_b V_1^*(b)
\end{equation*}
Therefore $a\mapsto \bK^{-1} \partial_a U_1^*(a)$ and $b\mapsto \bS^{-1} \partial_b V_1^*(b)$ are maximally monotone. Next, $\alpha \mapsto \partial_\alpha U_0^*(\alpha)$ and $\beta\mapsto \partial_\beta V_0^*(\beta)$ are maximally monotone by the convexity of $U_0^*,V_0^*$. Hence, $M$ is maximally monotone.


\end{document}